\documentclass[reqno,12pt]{amsart}
\usepackage[all]{xy}
\theoremstyle{plain}

\newtheorem{theorem}{Theorem}

\theoremstyle{definition}
\newtheorem{remark}[theorem]{Remark}

\newcommand{\C}{\mathbb{C}}

\newcommand{\F}{\widetilde{\mathcal F}}
\newcommand{\M}{\widetilde{M}}

\def\ka{\kappa}

\def\Dom{\operatorname{Dom}}

\def\sup{\operatorname{sup}}

\def\C{\mathbb C}
\def\R{\mathbb R}
\def\Z{\mathbb Z}

\def\N{\mathbb N}
\def\B{\mathcal B}

\def\A{{\mathcal A}} \def\B{{\mathcal B}}
 
\def\F{{\mathcal F}}

\def\tr{\operatorname{tr}}

\def\id{\operatorname{id}}
\def\H{\mathcal H}

\def\U{{\mathcal U}}
\def\<{\langle}
\def\>{\rangle}

\newcommand{\supp}{\operatorname{supp}}
\newcommand{\nc}{\newcommand}
\nc{\nt}{\newtheorem}
\nc{\gf}[2]{\genfrac{}{}{0pt}{}{#1}{#2}}
\nc{\mb}[1]{{\mbox{$ #1 $}}}
\nc{\real}{{\mathbb R}}
\nc{\comp}{{\mathbb C}}
\nc{\ints}{{\mathbb Z}}
\nc{\Ltoo}{\mb{L^2({\mathbf H})}}
\nc{\rtoo}{\mb{{\mathbf R}^2}}
\nc{\slr}{{\mathbf {SL}}(2,\real)}
\nc{\slz}{{\mathbf {SL}}(2,\ints)}
\nc{\su}{{\mathbf {SU}}(1,1)}
\nc{\so}{{\mathbf {SO}}}
\nc{\hyp}{{\mathbb H}}
\nc{\disc}{{\mathbf D}}
\nc{\torus}{{\mathbb T}}

\newcommand{\cH}{{\mathcal H}}
\nc{\ca}{{\mathcal A}}
\nc{\cag}{{{\mathcal A}^\Gamma}}
\nc{\cg}{{\mathcal G}}
\nc{\chh}{{\mathcal H}}
\nc{\ck}{{\mathcal B}}
\nc{\cm}{{\mathcal M}}
\nc{\cs}{{\mathcal S}}
\nc{\cz}{{\mathcal Z}}
\nc{\sind}{\sigma{\rm -ind}}

\begin{document}

\title[Spectral gaps for periodic Schr\"odinger operators]  {Spectral gaps for periodic
Schr\"odinger operators with strong magnetic fields}
\author{Yuri A. Kordyukov}
\address{Institute of Mathematics, Russian Academy of Sciences, Ufa,
Russia} \email{yuri@imat.rb.ru}

\thanks{}

\begin{abstract}
We consider Schr\"odinger operators $H^h = (ih\,d+{\bf A})^*
(ih\,d+{\bf A})$ with the periodic magnetic field ${\bf B}=d{\bf
A}$ on covering spaces of compact manifolds. Under some
assumptions on $\bf B$, we prove that there are arbitrarily large
number of gaps in the spectrum of these operators in the
semiclassical limit of strong magnetic field $h\to 0$.
\end{abstract}
\maketitle

\section*{Introduction}

Let $(M,g)$ be a closed Riemannian oriented manifold of dimension
$n\geq 2$, $\widetilde M$ be its universal cover and $\widetilde
g$ be the lift of $g$ to $\M$ so that $\widetilde g$ is a
$\Gamma$-invariant Riemannian metric on $\M$ where $\Gamma$
denotes the fundamental group of $M$ acting on $\M$ by the deck
transformations. Let $\bf B$ be a real-valued $\Gamma$-invariant
closed 2-form on $\M$. We assume that $\bf B$ is exact. Choose a
real-valued 1-form $\bf A$ on $\M$ such that $d{\bf A} = \bf B$.
Physically we can think of $\bf A$ as the electromagnetic vector
potential for a magnetic field $\bf B$.

We consider the magnetic Schr\"odinger operator $H_h$ given by
\[
H^h = (ih\,d+{\bf A})^* (ih\,d+{\bf A}),
\]
acting on the Hilbert space $\mathfrak{H}=L^2(\M)$. Here $h>0$ is
a semiclassical parameter, which is assumed to be small.

In local coordinates $X=(X_1,\ldots,X_n)$, we write the 1-form
$\bf A$ as
\[
{\bf A}= A_1(X)\,dX_1+\ldots+A_n(X)\,dX_n,
\]
the matrix of the Riemannian metric $\widetilde g$ as
$g(X)=(g_{jl}(X))_{1\leq j,l\leq n}$ and its inverse as
$(g^{jl}(X))_{1\leq j,l\leq n}$. If $|g(X)|=\det(g(X))$, then the
Schr\"odinger operator $H_h$ is given by
\[
H_h=\frac{1}{\sqrt{|g(X)|}}\sum_{1\leq j,l\leq n}\left(i h
\frac{\partial}{\partial X_j}+A_j(X)\right)\left[\sqrt{|g(X)|}
g^{jl}(X) \left(i h \frac{\partial}{\partial
X_l}+A_l(X)\right)\right].
\]

For any $x\in\M$ denote by $B(x)$ the anti-symmetric linear
operator on the tangent space $T_x\M$ associated with the 2-form
$\bf B$:
\[
{\widetilde g}_x(B(x)u,v)={\bf B}_x(u,v),\quad u,v\in T_x\M.
\]
The trace-norm $|B(x)|$ of $B(x)$ is given by the formula
\[
|B(x)|=[\operatorname{Tr}(B^*(x)\cdot B(x))]^{1/2}.
\]

We will assume that there exists an integer $k>0$ such that, if
$B(x_0)=0$ then
\[
C_1^{-1}d(x, x_0)^k\leq |B(x)| \leq C_1d(x, x_0)^k
\]
in some neighborhood of $x_0$ (here $d$ denotes the geodesic
distance on $\M$). We assume that there exists at least one zero
of $B$.

\begin{theorem}
\label{main0} Under current assumptions, there exists an
increasing sequence $\{\lambda_m, m\in\N \}$, satisfying
$\lambda_m\to\infty$ as $m\to\infty$, such that for any $a$ and
$b$, satisfying $\lambda_m<a<b < \lambda_{m+1}$ with some $m$, the
interval $[ah^{\frac{2k+2}{k+2}}, bh^{\frac{2k+2}{k+2}}]$ does not
meet the spectrum of $H^h$ for any $h>0$ small enough. It follows
that there exists arbitrarily large number of gaps in the spectrum
of $H^h$ provided the coupling constant $h$ is sufficiently small.
\end{theorem}

Here the sequence $\{\lambda_m\}$ appears as the set of
eigenvalues of a {\em model operator} $K^h$ associated to $H^h$.
This operator is defined as a direct sum of principal parts of
$H^h$ near the zeroes of $B$ in a fundamental domain (see
Section~\ref{s:model} for a precise definition). It is a
differential operator, which acts on the Hilbert space
$\mathfrak{H}_K= L^2({\mathbb R}^n)^N$ and has discrete spectrum
(here $n$ is the dimension of $M$ and $N$ denotes the number of
zeroes of $B$ that lie in a fundamental domain). Using a simple
scaling and gauge invariance, it can be shown that the operator
$K^h$ is unitarily equivalent to the operator
$h^{\frac{2k+2}{k+2}}K^1$. Therefore the operator
$h^{-\frac{2k+2}{k+2}}K^h$ has discrete spectrum independent of
$h$. This fact explains the appearance of a scaling factor
$h^{\frac{2k+2}{k+2}}$ in Theorem~\ref{main0}.

There exist a few examples of periodic partial differential
operators of the second order with spectral gaps (see, for
instance,
\cite{FK1,FK2,Fr,HempelHerbst95,HempelLienau00,HerbstNakamura99}
and a recent survey \cite{HempelPost02} and references therein).
In particular, in \cite{HempelHerbst95} Hempel and Herbst studied
magnetic Schr\"odinger operators
\[
H(\lambda \vec{a})=(-i\nabla -\lambda\vec{a}(x))^2
\]
in $L^2(\R^n)$, where $\vec{a}\in C^1(\R^n;\R^n)$ is a vector
potential and $\lambda\in\R$. Let $M=\{x\in\R^n:B(x)=0\}$, where
$B=d\vec{a}$ is the magnetic field associated with $\vec{a}$, and
$M_{\vec{a}}=\{x\in\R^n:\vec{a}(x)=0\}$. They proved that, if $B$
is periodic with respect to the lattice $\Z^n$, the set
$M\setminus M_{\vec{a}}$ has measure zero, the interior of $M$ is
non-empty and $M$ can be represented as $M=\cup_{j\in {\Z^n}} M_j$
(up to a set of measure zero) where the $M_j$ are pairwise
disjoint compact sets with $M_j=M_0+j$, then the spectrum of the
operator $H(\lambda \vec{a})$ has arbitrarily large number of gaps
provided the coupling constant $\lambda$ is sufficiently large.
The proof of this result is based on the fact that, as
$\lambda\to\infty$, $H(\lambda \vec{a})$ converges in norm
resolvent sense to the Dirichlet Laplacian $-\Delta_M$ on the
closed set $M$. Since norm resolvent convergence implies
convergence of spectra, we immediately obtain that, as
$\lambda\to\infty$, the spectrum of $H(\lambda \vec{a})$
concentrates around the eigenvalues of $-\Delta_M$ and gaps opens
up in the spectrum of $H(\lambda \vec{a})$.

On the other hand, Hempel and Herbst also proved in
\cite{HempelHerbst95} that, if $M_{\vec{a}}$ has measure zero,
then, as $\lambda\to\infty$, $H(\lambda \vec{a})$ converges in
strong resolvent sense to the zero operator in $L^2(\R^n)$. So, in
this case, their method to produce operators with spectral gaps
does not work.

In this paper, we consider a particular case when $M_{\vec{a}}$
has measure zero. More precisely, Theorem~\ref{main0} states that
if $M_{\vec{a}}$ has measure zero and the magnetic field has a
regular behaviour near its zeroes, we still can produce examples
of magnetic Schr\"odinger operators $H(\lambda \vec{a})$ with
arbitrarily large number of gaps in their spectra. The proof of
Theorem~\ref{main0} is based on Theorem~\ref{main1} below.

Recall that the magnetic Schr\"odinger operator $H^h$ commutes
with a projective $(\Gamma, \sigma)$-action of the fundamental
group $\Gamma$, where $\sigma$ is the multiplier or $U(1)$-valued
2-cocycle on $\Gamma$ defining this projective action.

Consider the reduced twisted group $C^*$-algebra
$C^*_r(\Gamma,\bar\sigma)$ of the group $\Gamma$. If $\cH$ is a
Hilbert space, then let ${\mathcal K}(\cH)$ denote the algebra of
compact operators in $\cH$, and $\mathcal{K}=
\mathcal{K}(\ell^2(\N))$, where $\N=\{1,2,3,\ldots\}$. Let
$E^h(\lambda) = {\chi}_{(-\infty, \lambda]}(H^h)$  and
$E^{0}(\lambda) = \chi_{(-\infty, \lambda]}(K^h)$ denote the
spectral projections. One can define actions of the $C^*$-algebra
$C^*_r(\Gamma, \bar\sigma)\otimes \mathcal K $ in the Hilbert
spaces $\mathfrak{H}$ and $\ell^2(\Gamma)\otimes \mathfrak{H}_K$.
It can be shown that $E(\lambda )$ and $\id\otimes E^{0}(\lambda)$
are in $C^*_r(\Gamma, \bar\sigma)\otimes \mathcal K $.

Recall that two projections $P$ and $Q$ in a unital $*$-algebra
$\mathcal A$ are said to be {\it Murray-von Neumann equivalent} if
there is an element $V\in \mathcal A$ such that $P= V^*V$ and $Q=
VV^*$.

\begin{theorem}\label{main1} Assume that $\lambda\in \R$ does not
coincide with $\lambda_k$ for any $k$. There exists a $(\Gamma,
\sigma)$-equivariant isometry $U:{\mathfrak H}\to
\ell^2(\Gamma)\otimes {\mathfrak H}_K$ and a constant $h_0
>0$ such that for all $h\in (0, h_0)$, the spectral projections $U
E(h^{\frac{2k+2}{k+2}}\lambda ) U^*$ and $\id\otimes
E^0(h^{\frac{2k+2}{k+2}}\lambda)$ are Murray-von Neumann
equivalent in $C^*_r(\Gamma, \bar\sigma)\otimes \mathcal
K({\mathfrak H}_K)$.
\end{theorem}

The proofs of Theorems~\ref{main0} and \ref{main1} are based on
abstract operator-theoretic results on equivalence of projections
and existence of spectral gaps proved in \cite{KMS}. These results
were applied in \cite{KMS} to prove existence of arbitrarily large
number of gaps in the spectrum of periodic magnetic Schr\"odinger
operators
\[
H_\mu=(i\,d+{\bf A})^* (i\,d+{\bf A})+\mu^{-2}V(x)
\]
on the universal covering $\M$ of a compact manifold $M$ in the
limit of the strong electric field ($\mu\to 0$), where ${\bf
B}=d{\bf A}$ is a $\Gamma$-invariant closed 2-form on $\M$ and
$V\geq 0$ is a $\Gamma$-invariant Morse potential (see also
\cite{MS}).

The another important results, which we use in construction of the
model operator and in the proof of Theorem~\ref{main1}, are
connected with the study of Schr\"odinger operators with magnetic
wells and were obtained by Hellfer and Mohamed (=Morame) in
(\cite{HM}, see also \cite{HM01,HM02,HM04} for further
developments).

The paper is organized as follows. Section~\ref{s:abstract-equiv}
contains some background results from \cite{KMS}. In
Section~\ref{s:model}, we describe a construction of the model
operator $K^h$. Section~\ref{operator} provides some necessary
information on magnetic translations and related operator
algebras. Finally, in Section~\ref{s:model1} we give the proofs of
the main results.

I am very thankful to Bernard Helffer for bringing these problems
to my attention and useful discussions.

\section{General results on equivalence of projections and
existence of spectral gaps}\label{s:abstract-equiv}

In this section  we recall general results on equivalence of
projections and existence of spectral gaps proved in \cite{KMS}.

Let $\mathfrak A$ be a $C^*$-algebra, $\chh$ a Hilbert space
equipped with a faithful $\ast$-representation of ${\mathfrak A}$,
$\pi: {\mathfrak A}\to \ck(\chh)$. For simplicity of notation, we
will often identify the algebra ${\mathfrak A}$ with its image
$\pi({\mathfrak A})$.

Consider Hilbert spaces ${\mathcal H}_1$ and ${\mathcal H}_2$
equipped with inner products $(\cdot,\cdot)_1$ and
$(\cdot,\cdot)_2$. Assume that there are given unitary operators
${\mathcal V}_1 : {\mathcal H}_1\to {\mathcal H}$ and ${\mathcal
V}_2 : {\mathcal H}_2\to {\mathcal H}$. Using the unitary
isomorphisms ${\mathcal V}_1$ and ${\mathcal V}_2$, we get
representations $\pi_1$ and $\pi_2$ of ${\mathfrak A}$ in
${\mathcal H}_1$ and ${\mathcal H}_2$ accordingly,
$\pi_l(a)={\mathcal V}^{-1}_l\circ\pi(a)\circ {\mathcal V}_l,
l=1,2, a\in {\mathfrak A}$.

Consider (unbounded) self-adjoint operators $A_1$ in ${\mathcal
H}_1$ and $A_2$ in ${\mathcal H}_2$ with the domains $\Dom (A_1)$
and $\Dom (A_2)$ respectively. We will assume that
\begin{itemize}
\item the operators $A_1$ and $A_2$ are semi-bounded from below:
\begin{gather}\label{e:5}
(A_1u,u)_1\geq \lambda_{01}\|u\|_1^2,\quad u\in \Dom (A_1), \\
\label{e:6} (A_2u,u)_2\geq \lambda_{02}\|u\|_2^2,\quad u\in \Dom
(A_2),
\end{gather}
with some $\lambda_{01},\lambda_{02}\leq 0$;
\item for any $t>0$, the operators $e^{-tA_l}, l=1,2,$ belong to
$\pi_l({\mathfrak A})$.
\end{itemize}

Let ${\mathcal H}_0$ be a Hilbert space, equipped with injective
bounded linear maps $i_1:{\mathcal H}_0\to {\mathcal H}_1$ and
$i_2:{\mathcal H}_0\to {\mathcal H}_2$. Assume that there are
given bounded linear maps $p_1:{\mathcal H}_1\to {\mathcal H}_0$
and $p_2:{\mathcal H}_2\to {\mathcal H}_0$ such that $p_1\circ
i_1=\id_{{\mathcal H}_0}$ and $p_2\circ i_2=\id_{{\mathcal H}_0}$.
The whole picture can be represented by the following diagram
(note that this diagram is not commutative).

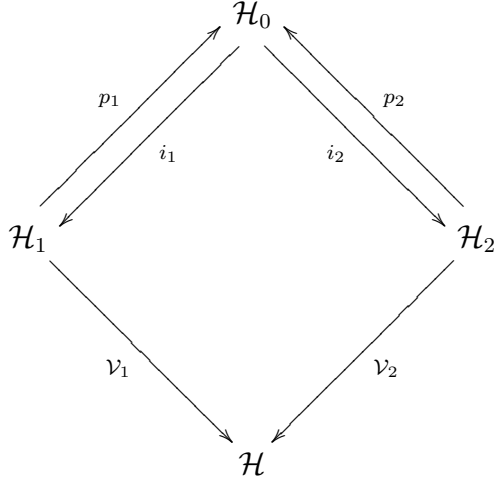
\begin{equation*}
\xymatrix @=8pc @ur { \cH_1 \ar@<1ex>[r]^{p_1}\ar[d]_{{\mathcal
V}_1} & \cH_0\ar@<-1ex>[d]_{i_2} \ar@<1ex>[l]^{i_1} \\ \cH & \cH_2
\ar@<-1ex>[u]_{p_2} \ar[l]^{{\mathcal V}_2}}
\end{equation*}

Consider a self-adjoint bounded operator $J$ in ${\mathcal H}_0$.
We assume that
\begin{itemize}
\item the operator ${\mathcal V}_2i_2Jp_1{\mathcal
V}^{-1}_1$ belongs to the von Neumann algebra $\pi({\mathfrak
A})''$;
\item $(i_2Jp_1)^*=i_1Jp_2$;
\item for any $a\in {\mathfrak A}$, the operator $\pi(a){\mathcal
V}_2(i_2Jp_1){\mathcal V}^{-1}_1$ belongs to $\pi({\mathfrak A})$.
\end{itemize}

Since the operators $i_l:{\mathcal H}_0\to {\mathcal H}_l, l=1,2,$
are bounded and have bounded left-inverse operators $p_l$, they
are topological monomorphisms, i.e. they have closed image and the
maps $i_l:{\mathcal H}_0\to {\rm Im}\, i_l$ are topological
isomorphisms. Therefore, we can assume that the estimate
\begin{equation}\label{e:rho}
\rho^{-1}\|i_2Ju\|_2\leq \|i_1Ju\|_1\leq \rho\|i_2Ju\|_2, \quad
u\in {\mathcal H}_0,
\end{equation}
holds with some $\rho>1$ (depending on $J$).

Define the bounded operators $J_l$ in ${\mathcal H}_l, l=1,2,$ by
the formula $J_l=i_lJp_l$. We assume that
\begin{itemize}
\item the operator $J_l, l=1,2,$ maps the domain of $A_l$ to itself;
\item $J_l$ is self-adjoint, and $0\leq J_l\leq \id_{{\mathcal H}_l}, l=1,2$;
\item for $u\in {\mathcal H}_0$, $i_1Ju\in \Dom(A_1)$ iff $i_2Ju\in \Dom(A_2)$.
\end{itemize}
Denote $D=\{u\in {\mathcal H}_0 : i_1Ju\in \Dom(A_1)\}=\{u\in
{\mathcal H}_0 : i_2Ju\in \Dom(A_2)\}.$

Introduce a self-adjoint positive bounded linear operator $J'_l$
in ${\mathcal H}_l$ by the formula
$J_l^2+{J'_l}{}^2=\id_{{\mathcal H}_l}$. We assume that
\begin{itemize}
\item the operator $J'_l, l=1,2,$ maps the domain of $A_l$ to itself;
\item the operators $[J_l,[J_l,A_l]]$ and $[J'_l,[J'_l,A_l]]$ extend to bounded
operators in ${\mathcal H}_l$, and
\begin{equation}\label{e:15}
\max (\|[J_l,[J_l,A_l]]\|_l,\, \|[J'_l,[J'_l,A_l]]\|_l)\leq
\gamma_l, \quad l=1,2.
\end{equation}
\end{itemize}

Finally, we assume that
\begin{equation}\label{e:14}
(A_lJ'_lu,J'_lu)_l\geq \alpha_l \|J'_lu\|_l^2,\quad u\in
\Dom(A_l), \quad l=1,2,
\end{equation}
for some $\alpha_l>0$, and
\begin{align}\label{e:16}
(A_2i_2Ju,i_2Ju)_2\leq \beta_1
(A_1i_1Ju,i_1Ju)_1+\varepsilon_1\|i_1Ju\|_1^2, \quad u\in D, \\
\label{e:A1A2} (A_1i_1Ju,i_1Ju)_1\leq \beta_2
(A_2i_2Ju,i_2Ju)_2+\varepsilon_2\|i_2Ju\|_2^2, \quad u\in D,
\end{align}
for some $\beta_1,\beta_2\geq 1$ and $\varepsilon_1,
\varepsilon_2>0$.

Denote by $E_l(\lambda), l=1,2$, the spectral projection of the
operator $A_l$, corresponding to the semi-axis
$(-\infty,\lambda]$. We assume that there exists a faithful,
normal, semi-finite trace $\tau$ on $\pi({\mathfrak A})''$ such
that, for any $t>0$, the operators ${\mathcal
V}_le^{-tA_l}{\mathcal V}^{-1}_l, l=1,2,$ belong to
$\pi({\mathfrak A})$ and have finite trace. By standard arguments,
it follows that ${\mathcal V}_lE_l(\lambda){\mathcal V}^{-1}_l\in
\pi ({\mathfrak A})''$, and $\tau({\mathcal
V}_lE_l(\lambda){\mathcal V}^{-1}_l)<\infty$ for any $\lambda,
l=1,2$.

\begin{theorem}\label{t:equivalence}
Under current assumptions, let $b_1>a_1$ and
\begin{align}\label{e:a2}
a_2&=\rho\left[ \beta_1 \left(a_1+\gamma_1+
\frac{(a_1+\gamma_1-\lambda_{01})^2}{\alpha_1-a_1-\gamma_1}\right)+
\varepsilon_1\right],\\ \label{e:b2}
b_2&=\frac{\beta_2^{-1}(b_1\rho^{-1}-\varepsilon_2)(\alpha_2-\gamma_2)
-\alpha_2\gamma_2+2\lambda_{02}\gamma_2-\lambda^2_{02}}
{\alpha_2-2\lambda_{02}+\beta_2^{-1}(b_1\rho^{-1}-\varepsilon_2)}.
\end{align}
Suppose that $\alpha_1>a_1+\gamma_1$, $\alpha_2>b_2+\gamma_2$ and
$b_2>a_2$. If the interval $(a_1,b_1)$ does not intersect with the
spectrum of $A_1$, then:

(1) the interval $(a_2,b_2)$ does not intersect with the spectrum
of $A_2$;

(2) for any $\lambda_1\in (a_1,b_1)$ and $\lambda_2\in (a_2,b_2)$,
the projections ${\mathcal V}_1E_1(\lambda_1){\mathcal V}^{-1}_1$
and ${\mathcal V}_2E_2(\lambda_2){\mathcal V}^{-1}_2$ belong to
${\mathfrak A}$ and are Murray-von Neumann equivalent in
${\mathfrak A}$.
\end{theorem}

\begin{remark}\label{b1b2}
Since $\rho>1, \beta_1\geq 1, \gamma_1>0$ and $\varepsilon_1>0$,
we, clearly, have $a_2>a_1$. The formula (\ref{e:b2}) is
equivalent to the formula
\[
b_1=\rho\left[ \beta_2 \left(b_2+\gamma_2+
\frac{(b_2+\gamma_2-\lambda_{02})^2}{\alpha_2-b_2-\gamma_2}\right)+
\varepsilon_2\right],
\]
which is obtained from (\ref{e:a2}), if we replace $\alpha_1,
\beta_1, \gamma_1, \varepsilon_1, \lambda_{01}$ by $\alpha_2,
\beta_2, \gamma_2, \varepsilon_2, \lambda_{02}$ accordingly and
$a_1$ and $a_2$ by $b_2$ and $b_1$ accordingly. In particular,
this implies that $b_1>b_2$.
\end{remark}

\section{The model operator}
\label{s:model} Here we will give a construction of the model
operator, using ideas of \cite{HM}. We will use notation of
Introduction. Choose a fundamental domain $\F\subset\widetilde M$
so that there is no zeros of $B$ on the boundary of $\F$. This is
equivalent to saying that  the translations
$\{\gamma\F,\;\gamma\in\Gamma\}$ cover the set of all zeros of
$B$. Let $\{\bar x_j|\,j=1,\dots,N\}$ denote all the zeros of $B$
in $\F$; $\bar x_i\ne\bar x_j$ if $i\ne j$.

The {\em model operator} $K^h$ associated with $H^h$ is an
operator in $L^2({\mathbb R}^n)^N$ given by $$ K^h = \oplus_{1\le
j\le N} K^h_j, $$ where $K^h_j$ is an unbounded self-adjoint
differential operator in $L^2({\mathbb R}^n)$ which corresponds to
the zero $\bar x_j$. Let us fix local coordinates $f_j: U(\bar
x_j)\to \R^n$ on $\M$ defined in a small neighborhood $U(\bar
x_j)$ of $\bar x_j$ for every $j=1,\dots,N$. We assume that
$f_j(\bar x_j)=0$ and the image $f_j(U(\bar x_j))$ is a fixed ball
$B=B(0,r)\subset \R^n$ centered at the origin $0$.

Write the 2-form $\bf B$ in the local coordinates as
\[
{\bf B}_j(X)=\sum_{1\leq l<m\leq n} b_{lm}(X)\,dX_l\wedge dX_m,
\quad X=(X_1,\ldots,X_n)\in B(0,r).
\]
The 1-form $\bf A$ is written in the local coordinates as a 1-form
${\bf A}_j$ on $B(0, r)$. By \cite{HM88}, there exists a real
function $\theta_j\in C^{\infty}(B(0, r))$ such that
\[
|{\bf A}_j(X)-d\theta_j(X)|\leq C |X|^{k+1}, \quad X\in B(0, r).
\]
Write the 1-form ${\bf A}_j-d\theta_j$ as
\[
{\bf A}_{j}(X)-d\theta_j(X)=\sum_{l=1}^n a_{l}(X)\,dX_l, \quad
X\in B(0, r).
\]
Let ${\bf A}_{1,j}$ be a 1-form on $\R^n$ with polynomial
coefficients given by
\[
{\bf A}_{1,j}(X)=\sum_{l=1}^n\sum_{|\alpha|=k+1}
\frac{X^\alpha}{\alpha !} \frac{\partial^\alpha a_{l}}{\partial
X^\alpha}(0)\,dX_l,  \quad X\in\R^n.
\]
So we have
\[
d{\bf A}_{1,j}(X) ={\bf B}^0_j(X), \quad X\in\R^n,
\]
where ${\bf B}^0_j$ is a closed 2-form on $\R^n$ with polynomial
coefficients defined by
\[
{\bf B}^0_j(X)=\sum_{1\leq l<m\leq
n}\sum_{|\alpha|=k}\frac{X^\alpha}{\alpha !}\frac{\partial^\alpha
b_{lm}}{\partial X^\alpha}(0)\,dX_l\wedge dX_m, \quad X\in\R^n.
\]
Take any extension of the function $\theta_j$ to a smooth,
compactly supported function in $\R^n$ denoted also by $\theta_j$
and put
\[
{\bf A}^0_j(X)= {\bf A}_{1,j}(X)+ d\theta_j(X), \quad X\in\R^n.
\]
Then we still have
\[
d{\bf A}^0_j(X) ={\bf B}^0_j(X), \quad X\in\R^n,
\]
and, moreover,
\begin{gather}
|{\bf B}_j(X)- {\bf B}^0_j(X)|\leq C|X|^{k+1}, \quad X\in B(0,
r),\\ |{\bf A}_{j}(X)- {\bf A}^0_j(X)|\leq C|X|^{k+2}, \quad X\in
B(0, r). \label{e:a}
\end{gather}

Then $K^h_j$ is the self-adjoint differential operator with
asymptotically polynomial coefficients in $L^2(\R^n)$  given by
\[
K_j^h =  (i h\,d+{\bf A}^0_j)^* (i h\,d+{\bf A}^0_j),
\]
where the adjoint is taken with respect to a Hilbert structure in
$L^2(\R^n)$ given by the flat Riemannian metric $(g_{lm}(0))$ in
$\R^n$. If we write ${\bf A}^0_j$ as
\[
{\bf A}^0_j=A^0_{j,1}\, dX_1+\ldots+ A^0_{j,n}\,dX_n,
\]
then $K^h_j$ is given as
\[
K_j^h=\sum_{1\leq l,m\leq n} g^{lm}(0) \left(i h
\frac{\partial}{\partial X_l}+A^0_{j,l}(X)\right)\left(i h
\frac{\partial}{\partial X_m}+A^0_{j,m}(X)\right).
\]
The operator $K^h_j$ has discrete spectrum (cf., for instance,
\cite{HelNo85,HM88}). By gauge invariance, the operator $K^h_j$ is
unitarily equivalent to the Schr\"odinger operator
\[
H_j^h = (ih\,d+{\bf A}_{1,j})^* (ih\,d+{\bf A}_{1,j}),
\]
associated with the homogeneous $1$-form ${\bf A}_{1,j}$. Using a
simple scaling $X\mapsto h^{\frac{1}{k+2}}X$, it can be shown that
the operator $H_j^h$ is unitarily equivalent to the operator
$h^{\frac{2k+2}{k+2}}H_j^1$. So we conclude that the operator
$h^{-\frac{2k+2}{k+2}}K^h$ has discrete spectrum independent of
$h$, which is denoted by $\{\lambda_m:m\in \N\}$,
$\lambda_1<\lambda_2<\lambda_2<\ldots$ (not taking into account
multiplicities). As it will be shown in Section~\ref{s:model1},
the sequence $\{\lambda_m:m\in \N\}$ is precisely what we need for
the proof of Theorem~\ref{main0}.

\section{Magnetic translations and related operator
algebras}\label{operator} In this section, we collect some
necessary facts on magnetic translations and related operator
algebras (see, for instance, \cite{MS,KMS} and references therein
for more details). As above, let $M$ be a compact connected
Riemannian manifold, $\Gamma$ be its fundamental group and
$p:{\widetilde M}\to M$ be its universal cover. Let ${\bf B}$ be a
closed $\Gamma$-invariant real-valued 2-form on $M$. Assume that
$\bf B$ is \emph{exact}. So ${\bf B}=d{\bf A}$ where ${\bf A}$ is
a 1-form on $\widetilde M$. We will assume without loss of
generality that ${\bf A}$ is real-valued.

The Hermitian connection ${\bf A}$ defines a projective $(\Gamma,
{\sigma})$-unitary representation on $L^2(\M)$, that is, the map
$T : \Gamma\to {\U} (L^2(\M))$, $\gamma\mapsto T_\gamma$, where
for any Hilbert space $\H$ we denote by $\U(\H)$ the group of all
unitary operators in $\H$, satisfying $$ T_e={\id}, \quad
T_{\gamma_1} T_{\gamma_2} = {\sigma}(\gamma_1,\gamma_2)
T_{\gamma_1 \gamma_2}, \quad \gamma_1, \gamma_2 \in \Gamma.$$ Here
$\sigma$ is a {\em multiplier} on $\Gamma$ i.e.
$\sigma:\Gamma\times\Gamma\to U(1)$ satisfies
\begin{itemize}
\item ${\sigma}(\gamma,e) = {\sigma}(e,\gamma)=1,\quad
\gamma\in\Gamma$;

\item ${\sigma}(\gamma_1,\gamma_2)
{\sigma}(\gamma_1\gamma_2, \gamma_3)=
{\sigma}(\gamma_1,\gamma_2\gamma_3)
{\sigma}(\gamma_2,\gamma_3),\quad \gamma_1, \gamma_2, \gamma_3\in
\Gamma$ \quad ({\em the cocycle relation}).
\end{itemize}

In other words one says that the map $\gamma\mapsto T_\gamma$
defines a $(\Gamma,\sigma)$-{\em action} in~$\H$. The operators
$T_\gamma$ are also called {\em magnetic translations}.

Denote by $\ell^2(\Gamma)$ the standard Hilbert space of
complex-valued $L^2$-functions on the discrete group $\Gamma$. For
any $\gamma\in \Gamma$, define a bounded operator $T_\gamma^L$ in
$\ell^2(\Gamma)$ by $$ T_\gamma^L f(\gamma') =
f(\gamma^{-1}\gamma') \bar\sigma(\gamma, \gamma^{-1}\gamma'),
\quad \gamma' \in \Gamma, \quad f\in \ell^2(\Gamma). $$ It is easy
to see that $$ T_e^L={\id}, \quad T_{\gamma_1}^L T_{\gamma_2}^L=
\bar\sigma(\gamma_1,\gamma_2)T_{\gamma_1 \gamma_2}^L, \quad
\gamma_1,\gamma_2\in\Gamma. $$ Also $$
(T_\gamma^L)^*=\sigma(\gamma, \gamma^{-1})T_{\gamma^{-1}}^L. $$
This means that $T_\gamma$ is a left $(\Gamma, \bar\sigma)$-action
on $\ell^2(\Gamma)$ (or, equivalently, a $(\Gamma,
\bar\sigma)$-unitary representation in $\ell^2(\Gamma)$).

Define a twisted group algebra $\C(\Gamma,\bar\sigma)$ which
consists of complex valued functions with finite support on
$\Gamma$, with the twisted convolution operation $$
(f*g)(\gamma)=\sum_{\gamma_1,\gamma_2:\gamma_1\gamma_2=
\gamma}f(\gamma_1)g(\gamma_2)\bar\sigma(\gamma_1,\gamma_2), $$ and
with the involution $$f^*(\gamma)=\sigma(\gamma,
\gamma^{-1})\overline{f({\gamma^{-1}})}. $$ Associativity of the
multiplication is equivalent to the cocycle condition. The basis
of $\C(\Gamma,\bar\sigma)$ as a vector space is formed by
$\delta$-functions $\{\delta_\gamma\}_{\gamma\in\Gamma}$,
$\delta_\gamma(\gamma')=1$ if $\gamma=\gamma'$ and $0$ otherwise.
We have $$ \delta_{\gamma_1} *
\delta_{\gamma_2}=\bar\sigma(\gamma_1,\gamma_2)\delta_{\gamma_1\gamma_2}.
$$

Note also that the $\delta$-functions
$\{\delta_\gamma\}_{\gamma\in\Gamma}$ form an orthonormal basis in
$\ell^2(\Gamma)$. It is easy to check that $$
T_\gamma^L\delta_{\gamma'}=\delta_{\gamma} *
\delta_{\gamma'}=\bar\sigma(\gamma,\gamma')\delta_{\gamma\gamma'}.
$$

It is clear that the correspondence $f\in \C(\Gamma,\bar\sigma)
\mapsto T^L(f)\in \B(\ell^2(\Gamma)),$ where $T^L(f)u =f*u, u\in
\ell^2(\Gamma),$ defines a $\ast$-representation of the twisted
group algebra $\C(\Gamma,\bar\sigma)$ in $\ell^2(\Gamma)$. The
weak closure of the image of $\C(\Gamma,\bar\sigma)$ in this
representation coincides with the {\em (left) twisted group von
Neumann algebra} $\A^L(\Gamma,\bar\sigma)$. The corresponding norm
closure is so called {\em reduced twisted group} $C^*$-{\em
algebra} which is denoted $C^*_r(\Gamma,\bar\sigma)$.

The von Neumann algebra $\A^L(\Gamma,\bar\sigma)$ can be described
in terms of the matrix elements. For any $A\in\B(\ell^2(\Gamma))$
denote $A_{x,y}=(A\delta_y,\delta_x), x,y\in\Gamma$ (which is a
matrix element of $A$). Then repeating standard arguments (given
in a similar situation e.g.\ in \cite{Shubin}) we can prove that
for any $A\in\B(\ell^2(\Gamma))$ the inclusion $A\in
\A^L(\Gamma,\bar\sigma)$ is equivalent to the relations $$
A_{x\gamma,y\gamma}=\bar\sigma(x,\gamma)\sigma(y,\gamma)A_{x,y}\;,\quad
x,y,\gamma\in\Gamma. $$ In particular, for any $A\in
\A^L(\Gamma,\bar\sigma)$, we have $$ A_{x\gamma ,x\gamma
}=A_{x,x}\;, \quad x,\gamma\in\Gamma. $$

A finite {\em von Neumann trace}
$\tr_{\Gamma,\bar\sigma}:\A^L(\Gamma,\bar\sigma)\to\C$ is defined
by the formula $$ \tr_{\Gamma,\bar\sigma} A=(A\delta_e,\delta_e).
$$ We can also write $\tr_{\Gamma,\bar\sigma} A=A_{\gamma,\gamma}=
\left(A\delta_\gamma,\delta_\gamma\right)$  for any $\gamma\in
\Gamma$ because the right hand side does not depend of $\gamma$.

\section{Proof of main results}\label{s:model1}
For the proof of the main theorem, we apply
Theorem~\ref{t:equivalence} in a following particular setting.
Take the $C^*$ algebra ${\mathfrak A}$ to be $C^*_r(\Gamma,
\bar\sigma)\otimes {\mathcal K}$. Let $\mathcal H$ be the Hilbert
space $\ell^2(\Gamma)\otimes \ell^2(\N)$. Put ${\mathcal
H}_1=\ell^2(\Gamma)\otimes L^2({\mathbb R}^n)^N$ and ${\mathcal
H}_2=L^2(\M)$. Choose an arbitrary unitary isomorphism
$V_1:L^2({\mathbb R}^n)^N\to \ell^2(\N)$ and define an unitary
operator ${\mathcal V}_1 : {\mathcal H}_1\to {\mathcal H}$ as
${\mathcal V}_1=\id\otimes V_1$.

As in \cite{KMS}, define a $(\Gamma, \sigma)$-equivariant isometry
${\bf U} : L^2(\M)\cong \ell^2(\Gamma)\otimes L^2(\mathcal{F})$ by
the formula $$ {\bf U} (\phi) = \sum_{\gamma\in
\Gamma}\delta_\gamma\otimes i^*(T_\gamma \phi), \qquad  \phi \in
L^2(\M), $$ where $i: \mathcal{F} \to \M$ denotes the inclusion
map. Choose an arbitrary unitary isomorphism $V_2:L^2(\F)\to
\ell^2(\N)$. Then a unitary operator ${\mathcal V}_2 : {\mathcal
H}_2\to {\mathcal H}$ is defined as ${\mathcal V}_2=(\id\otimes
V_2)\circ {\bf U}$.

Let $\pi$ be the representation of the algebra ${\mathfrak A}$ in
$\mathcal H$ given by the tensor product of the representation
$T^L$ of $C^*_r(\Gamma, \bar\sigma)$ on $\ell^2(\Gamma)$ and the
standard representation of ${\mathcal K}$ in $\ell^2(\N)$. So we
have $\pi(C^*_r(\Gamma, \bar\sigma)\otimes {\mathcal K})\subset
{\mathcal A}^L(\Gamma,\bar\sigma) \otimes \B(\ell^2(\N))$ and
$\pi(C^*_r(\Gamma, \bar\sigma)\otimes {\mathcal K})''\cong
{\mathcal A}^L(\Gamma,\bar\sigma) \otimes \B(\ell^2(\N))$. Using
the unitary isomorphisms ${\mathcal V}_1$ and ${\mathcal V}_2$, we
get representations $\pi_1$ and $\pi_2$ of ${\mathfrak A}$ in
${\mathcal H}_1$ and ${\mathcal H}_2$ accordingly,
$\pi_l(a)={\mathcal V}^{-1}_l\circ\pi(a)\circ {\mathcal V}_l,
l=1,2, a\in {\mathfrak A}$.

Define a trace $\tau$ on ${\mathcal A}^L(\Gamma,\bar\sigma)
\otimes \B(\ell^2(\N))$ as the tensor product of the finite von
Neumann trace $\tr_{\Gamma,\bar\sigma}$ on $\A^L(\Gamma,
\bar\sigma)$ and the standard trace on $\B(\ell^2(\N))$.

Consider self-adjoint, semi-bounded from below operators $A_1$ in
${\mathcal H}_1$ and $A_2$ in ${\mathcal H}_2$:
\[
A_1=\id \otimes h^{-\frac{2k+2}{k+2}} K^h,\quad
A_2=h^{-\frac{2k+2}{k+2}} H^h.
\]
Clearly, we have
\[e^{-tA_1}=\id \otimes
e^{-th^{-\frac{2k+2}{k+2}}K^h}\in\pi_1({\mathfrak A})\cong
C^*_r(\Gamma, \bar\sigma)\otimes {\mathcal K}(\mathfrak{H}_K)
\]
with $\tau(e^{-tA_1})<\infty$ for any $t>0$. As shown in
\cite{KMS}, for any $t>0$, the operator $e^{-tA_2}$ belongs to
$\pi_2({\mathfrak A})$ and $\tau(e^{-tA_2})<\infty$. Remark that,
in notation of Theorem~\ref{main1},
\[
E_1(\lambda)=\id\otimes
E^0(h^{\frac{2k+2}{k+2}}\lambda), \quad
E_2(\lambda)=E(h^{\frac{2k+2}{k+2}}\lambda).
\]

We will use notation of Section~\ref{s:model}. Let
\[
{\mathcal H}_0=\ell^2(\Gamma)\otimes \left(\oplus_{j=1}^N
L^2(U(\bar x_j))\right).
\]
An inclusion $ i_1: {\mathcal H}_0 \to {\mathcal H}_1$ is defined
as $i_1=\id \otimes j_1$, where $j_1$ is the inclusion
\[
\oplus_{j=1}^NL^2(U(\bar x_j))\cong L^2(B(0, r))^N \hookrightarrow
L^2({\mathbb R}^n)^N
\]
given by the chosen local coordinates. An inclusion $i_2:
{\mathcal H}_0 \to {\mathcal H}_2$ is defined as $i_2={\bf
U}^{*}\circ (\id\otimes j_2)$, where $j_2$ is the natural
inclusion
\[
\oplus_{j=1}^NL^2(U(\bar x_j))\hookrightarrow L^2(\F).
\]

The operator $p_1:{\mathcal H}_1\to {\mathcal H}_0$ is defined as
$p_1=\id\otimes r_1$, where $r_1$ is the restriction operator
\[
L^2({\mathbb R}^n)^N\to L^2(B(0, r))^N \cong
\oplus_{j=1}^NL^2(U(\bar x_j)).
\]
The operator $p_2:{\mathcal H}_1\to {\mathcal H}_0$ is defined as
$p_2=(\id\otimes r_2)\circ {\bf U}$, where $r_2: L^2(\F)\to
\oplus_{j=1}^NL^2(U(\bar x_j))$ is the restriction operator.

Fix a function $\phi\in C_c^\infty(\R^n)$ such that $0\leq\phi\leq
1$, $\phi(x)=1$ if $|x|\leq 1$, $\phi(x)=0$ if $|x|\geq 2$, and
$\phi'=(1-\phi^2)^{1/2}\in C^\infty(\R^n)$. Fix a number $\ka>0$,
which we shall choose later. For any $h>0$ define
$\phi^{(h)}(x)=\phi(h^{-\ka}x)$. For any $h>0$ small enough, let
$\phi_j=\phi^{(h)}\in C^\infty_c(U(\bar x_j))$ in the fixed
coordinates near $\bar x_j$. Denote also $\phi_{j,\gamma} =
(\gamma^{-1})^*\phi_j$. (This function is supported in the
neighborhood $U(\gamma\bar x_j)=\gamma(U(\bar x_j))$ of
$\gamma\bar x_j$.) We will always take $h\in (0,h_0)$ where $h_0$
is sufficiently small, so in particular the supports of all
functions $\phi_{j,\gamma}$ are disjoint.

Let $\Phi\in  C^\infty(\bigcup_{j=1}^N U(\bar x_j))$ be equal to
$\phi_j$ on $U(\bar x_j)$, $j=1,2,\ldots,N$. Consider a $(\Gamma,
\sigma)$-equivariant, self-adjoint, bounded operator $J$ in
${\mathcal H}_0$ defined as $J=\id\otimes \Phi$, where $\Phi$
denotes the multiplication operator by the function $\Phi$ in the
space $\oplus_{j=1}^NL^2(U(\bar x_j))$.

It is clear that
\[
{\mathcal V}_2i_2Jp_1{\mathcal V}_1^{-1}=\id\otimes V_2j_2\Phi
r_1V_1^{-1},
\]
and $j_2\Phi r_1$ is a bounded operator from $L^2({\mathbb
R}^n)^N$ to $L^2(\F)$ given as the composition
\[
L^2({\mathbb R}^n)^N\to L^2(B(0, r))^N \cong
\oplus_{j=1}^NL^2(U(\bar x_j)) \stackrel{\Phi}\longrightarrow
\oplus_{j=1}^NL^2(U(\bar x_j)) \hookrightarrow L^2(\F).
\]
Hence, the operator ${\mathcal V}_2i_2Jp_1{\mathcal V}^{-1}_1$
belongs to the von Neumann algebra $\pi({\mathfrak A})'' \cong
{\mathcal A}^L (\Gamma,\bar\sigma)\otimes {\mathcal
B}(\ell^2(\N))$, and, for any $a\in {\mathfrak A}$, the operator
$\pi(a){\mathcal V}_2(i_2Jp_1){\mathcal V}^{-1}_1$ belongs to
$\pi({\mathfrak A})$.

Similarly, we have
\[
i_1Jp_2=\id\otimes j_1\Phi r_2,
\]
and $j_1\Phi r_2$ is a bounded operator from $L^2(\F)$ to
$L^2({\mathbb R}^n)^N$ given as the composition
\[
L^2(\F)\to \oplus_{j=1}^NL^2(U(\bar
x_j))\stackrel{\Phi}\longrightarrow \oplus_{j=1}^NL^2(U(\bar x_j))
\cong L^2(B(0, r))^N \hookrightarrow L^2({\mathbb R}^n)^N.
\]
So we have $(i_2Jp_1)^*=i_1Jp_2$.

We will use local coordinates near $\bar{x}_j$ such that the
Riemannian volume element at the point $\bar{x}_j$ coincides with
the Euclidean volume element given by the chosen local
coordinates. Then the estimate \eqref{e:rho} holds with
\begin{equation}\label{e:rho1}
\rho=1+O(h^\kappa).
\end{equation}

Denote by the same letters $\phi$ and $\phi'$ the multiplication
operators in $L^2({\mathbb R}^n)$ by the functions $\phi$ and
$\phi'$ accordingly. Let $\Phi_1$ and $\Phi'_1$ be the bounded
operators in $L^2({\mathbb R}^n)^N\cong L^2({\mathbb R}^n)\otimes
\C^N$ given by $\Phi_1=\phi\otimes \id_{\C^N}$ and
$\Phi'_1=\phi'\otimes \id_{\C^N}$. Then we have $J_1=\id\otimes
\Phi_1 $ and $J'_1=\id\otimes \Phi'_1$ in $\ell^2(\Gamma)\otimes
L^2({\mathbb R}^n)^N$.

Let $\Phi_\gamma \in  C^\infty(\widetilde M)$ be equal to
$\phi_{j,\gamma}$ on $U(\gamma\bar x_j)$, $j=1,2,\ldots,N$, and
$0$ otherwise. Put $\Phi_2=\sum_{\gamma\in\Gamma}\Phi_\gamma \in
C^{\infty}(\widetilde M)$. Define a function $\Phi'_2\in
C^{\infty}(\widetilde M), \Phi'_2\geq 0$ by the equation
$(\Phi_2)^2+(\Phi'_2)^2=1\ \text{in}\ C^{\infty}(\widetilde M)$.
The operators $J_2$ and $J'_2$ are given by the multiplication
operators by the functions $\Phi_2$ and $\Phi'_2$ in
$L^2(\widetilde{M})$ respectively.

The estimate~\eqref{e:15} hold with
\begin{equation}\label{e:gamma1}
\gamma_l =O(h^{-\frac{2k+2}{k+2}+2-2\kappa}),\quad l=1,2.
\end{equation}
Indeed, for any $j=1,2,\ldots,N$, the principal symbol
$a_{1,j}^{(2)}\in C^\infty(T^*{\mathbb R}^n)$ of $K^h_j$ is given
by
\begin{equation*}
a_{1,j}^{(2)}(x,\xi)=h^2\sum_{i,k=1}^n g^{ik}(\bar x_j)\xi_i\xi_k,
\quad (x,\xi)\in T^*{\mathbb R}^n.
\end{equation*}
So we have
\[ [J_1,[J_1,A_1]]=-h^{-\frac{2k+2}{k+2}+2}\id \otimes \left( \oplus_{1\leq j\leq
N} a_{1,j}^{(2)}(x,d\phi(x))\right)\] and \[[J'_1,[J'_1,A_1]]=
-h^{-\frac{2k+2}{k+2}+2}\id \otimes \left(\oplus_{1\leq j\leq N}
a_{1,j}^{(2)}(x,d\phi'(x))\right)\] in $\ell^2(\Gamma)\otimes
L^2({\mathbb R}^n)^N$. Similarly, the principal symbol
$a_2^{(2)}\in C^\infty(T^*\widetilde M)$ of $H^h$ is given by
\[
a_{2}^{(2)}(x,\xi)=h^2\sum_{i,k=1}^n g^{ik}(x)\xi_i\xi_k, \quad
(x,\xi)\in T^*\widetilde M.
\]
So the operators $[J_2,[J_2,A_2]], [J'_2,[J'_2,A_2]]$ are the
multiplication operators in $L^2(\widetilde M)$ by the functions
$- h^{-\frac{2k+2}{k+2}+2} a_2^{(2)}(x,d\Phi_2(x))$ and
$-h^{-\frac{2k+2}{k+2}+2} a_2^{(2)}(x,d\Phi'_2(x))$ accordingly.
Therefore,
\begin{align*}
\gamma_1&=h^{-\frac{2k+2}{k+2}+2} \max_{j=1,2,\ldots,N} \max
\left( \sup_{ x\in {\mathbb R}^n} (a^{(2)}_{1,j}(x,d\phi(x))),
\sup_{x\in {\mathbb R}^n} (a^{(2)}_{1,j}(x,d\phi'(x)))\right)\\ &
=O(h^{-\frac{2k+2}{k+2}+2-2\kappa}),\\
\gamma_2&=h^{-\frac{2k+2}{k+2}+2} \max \left(\sup_{x\in \M}
(a^{(2)}_2(x,d\Phi_2(x))), \sup_{x\in \M}
(a^{(2)}_2(x,d\Phi'_2(x)))\right)\\ &
=O(h^{-\frac{2k+2}{k+2}+2-2\kappa}).
\end{align*}
The estimates \eqref{e:14} hold with
\begin{equation}\label{e:alpha1}
\alpha_l=O(h^{-\frac{2k+2}{k+2}+k\kappa+1}),\quad l=1,2.
\end{equation}
Indeed, let $q_{0,j}^h$ denote the quadratic hermitian form
associated to $K_j^h$,
\[
q_{0,j}^h(u)=(K_j^hu,u)=\int_{{\mathbb R}^n} |ih\,d u+{\bf
A}_j^0u|^2\sqrt{g(0)}\,dx.
\]
Consider the operator $H_j^h =  (ih\,d+{\bf A}_{1,j})^*
(ih\,d+{\bf A}_{1,j})$ and denote by $q_{1,j}^h$ the quadratic
hermitian form associated to this operator:
\[
q_{1,j}^h(u)=(H_j^hu,u)=\int_{{\mathbb R}^n} |ih\,d u+{\bf
A}_{1,j}u|^2\sqrt{g(0)}\,dx,
\]
By \cite[Theorem 4.4]{HM}, there exists a constant $C_j>0$ such
that for any $h>0$,
\begin{equation}\label{e:99}
h\int_{\R^n} |B_{0,j}(x)|\, |u(x)|^2\,dx \leq C_j q_{1,j}^h(u),
\quad u\in C^\infty_c(\R^n).
\end{equation}
Using gauge invariance and (\ref{e:99}), we get
\begin{equation}\label{e:100}
q_{0,j}^h(u)  = q_{1,j}^h(e^{i\theta_j}u)
 \geq \frac{h}{C_j}\int_{\R^n} |B_{0,j}(x)|\, |e^{i\theta_j} u(x)|^2\,dx
 = \frac{h}{C_j}\int_{\R^n} |B_{0,j}(x)|\, |u(x)|^2\,dx .
\end{equation}
Similarly, let $q^h$ be the quadratic hermitian form associated to
$H^h$,
\[
q^h(u)=(H^hu,u)=\int_{\widetilde M} |ih\,d u+{\bf A}u|^2\,d\mu(x),
\]
where $d\mu$ denotes the Riemannian volume form on $\M$. By an
easy modification of the proof of Theorem 4.5 in \cite{HM}, one
can show that there exists a constant $C_0>0$, and for any $
\varepsilon\in (0,1)$, there exists a constant $C_\varepsilon>0$
such that for any $h\in (0,h_0]$,
\begin{equation}\label{e:107}
h\int_{\widetilde M} |B(x)|\, |u(x)|^2\,d\mu(x) \leq
C_0\left[q^h(u)+C_\varepsilon h^{2-\varepsilon}\|u\|^2\right],
\quad u\in C^\infty_c(\M).
\end{equation}
Assume that $k\kappa<2$ and take $\varepsilon\in (0,1)$ so that
$k\kappa<2-\varepsilon$. Since $|B(x)|\geq Ch^{k\kappa}$ for $x\in
\supp \Phi'_l, l=1,2$, the estimates \eqref{e:100} and
\eqref{e:107} easily imply the desired estimates \eqref{e:14}.

The constants $\lambda_{0l}, l=1,2,$ can be chosen to be
independent of $h$. Namely, one can take
\begin{equation}\label{e:lambda1}
\lambda_{01}=\lambda_{02}=0.
\end{equation}

Finally, the estimates \eqref{e:16} and \eqref{e:A1A2} hold with
\begin{equation}\label{e:beta1}
\beta_l=1+O(h^\kappa), \quad \varepsilon_l=O(h^{2\kappa
(k+2)-\kappa-\frac{2k+2}{k+2}}), \quad l=1,2.
\end{equation}
Using the inequality $|a+b|^2\leq |a|^2+2|a| |b|+|b|^2\leq
(1+\varepsilon )|a|^2+(1+\varepsilon^{-1} )|b|^2$ with
$\varepsilon=h^\kappa$, we get
\begin{align*}
q^h(\phi u)=&(H^h(\phi u),\phi u)=\int_{U(\bar x_j)} |ih\,d (\phi
u)+{\bf A}\phi u|^2\,d\mu(x) \\ \leq &(1+h^\kappa)\int_{B(0,r)}
|ih\,d (\phi u)+{\bf A}_j\phi u|^2\,\sqrt{g(0)}\,dx \\ \leq &
(1+h^\kappa)\int_{B(0,r)} |ih\,d (\phi u)+{\bf A}_j^0\phi
u|^2\,\sqrt{g(0)}\,dx \\ & + ch^{-\kappa} \int_{B(0,r)} |({\bf
A}_j-{\bf A}_j^0)\phi u|^2\,\sqrt{g(0)}\,dx.
\end{align*}

By (\ref{e:a}), we have
\[
\int_{B(0,r)} |({\bf A}_j-{\bf A}_j^0)\phi u|^2\,\sqrt{g(0)}\,dx
\leq C h^{2\kappa (k+2)} \int_{B(0,r)} |\phi u|^2 \, \sqrt{g(0)}
\, dx,
\]
that completes the proof of \eqref{e:16}. The proof of
\eqref{e:A1A2} is similar.

Now we complete the proofs of Theorems~\ref{main0}
and~\ref{main1}. As above, take $\{\lambda_m:m\in \N\}$,
$\lambda_1<\lambda_2<\lambda_2<\ldots,$  to be the spectrum of
the operator $h^{-\frac{2k+2}{k+2}}K^h$ (without taking into
account multiplicities), which is independent of $h$. Take any
$a$ and $b$ such that $\lambda_m<a<b < \lambda_{m+1}$ with some
$m$. Clearly, the spectrum of the operator $A_1$ coincides with
the spectrum of the operator $h^{-\frac{2k+2}{k+2}}K^h$.
Therefore, the interval $[a,b]$ does not intersect with the
spectrum of $A_1$. Take any open interval $(a_1,b_1)$ that
contains $[a,b]$ and does not intersect with the spectrum of
$A_1$. Using the estimates \eqref{e:rho1}, \eqref{e:gamma1},
\eqref{e:alpha1}, \eqref{e:lambda1} and \eqref{e:beta1}, one can
see that, for $a_2$ and $b_2$ given by \eqref{e:a2} and
\eqref{e:b2}, we have
\begin{equation}\label{e:estimate}
a_2=a_1+O(h^{s}),\quad b_2=b_1+O(h^{s}), \quad h\to 0,
\end{equation}
where
$s=\min\{(2k+3)\kappa-\frac{2k+2}{k+2},-\frac{2k+2}{k+2}+2-2\kappa
\}$. The best possible value of $s$ which is
\[s=\max_\kappa\min \{(2k+3)\kappa-\frac{2k+2}{k+2},-\frac{2k+2}{k+2}+2-2\kappa
\}=\frac{2}{(2k+5)(k+2)}\] is attained when
$\kappa=\frac{2}{2k+5}$.

Hence, if $h>0$ is small enough, we have $\alpha_1>a_1+\gamma_1$,
$\alpha_2>b_2+\gamma_2$, $b_2>a_2$ and the interval $(a_2,b_2)$
contains $[a,b]$. By Theorem~\ref{t:equivalence}, we conclude that
the interval $(a_2,b_2)$ does not intersect with the spectrum of
$A_2$, that completes the proof of Theorem~\ref{main0}. Moreover,
we have that, for any $\lambda_1\in (a_1,b_1)$ and $\lambda_2\in
(a_2,b_2)$, the spectral projections ${\mathcal
V}_1E_1(\lambda_1){\mathcal V}^{-1}_1$ and ${\mathcal
V}_2E_2(\lambda_2){\mathcal V}^{-1}_2$ are equivalent in
${\mathfrak A}$. Putting $U={\mathcal V}^{-1}_1{\mathcal V}_2$, we
get the desired Murray - von Neumann equivalence of
$E_1(\lambda_1)=\id\otimes E^0(\lambda)$ and ${\mathcal
V}^{-1}_1{\mathcal V}_2E_2(\lambda_2){\mathcal V}^{-1}_2{\mathcal
V}_1=U E(\lambda) U^{-1}$ in $\pi_1(\mathfrak A)= C^*_r(\Gamma,
\bar\sigma)\otimes \mathcal K(L^2({\mathbb R}^n)^N)$.

\end{document}